\documentclass[12pt]{amsart}
\usepackage{amsmath,bm}
\usepackage{amsfonts}
\usepackage{amssymb}
\usepackage{amsthm}
\usepackage{amscd}
\usepackage{a4wide}

\usepackage{amsmath}
\usepackage{amssymb}
\usepackage{amsthm}
\usepackage{graphics}
\usepackage{eucal}
\usepackage{mathrsfs}

\usepackage[utf8]{inputenc} 
\usepackage{textcomp} 
\usepackage{graphicx}  
\usepackage{epstopdf} 
\usepackage{flafter}  
\usepackage{topcapt}   
\usepackage{bm}  
\usepackage{tikz}
\usepackage[pdftex,bookmarks,colorlinks,breaklinks]{hyperref}  
\usepackage{color}
\usepackage{memhfixc}  
\usepackage{pdfsync}  
\usepackage{amsmath,amsthm,indentfirst,marvosym}
\usepackage{amssymb}
\usepackage{amsfonts}
\theoremstyle{plain}
\newtheorem*{maintheorem*}{Main Theorem}
\newtheorem*{thm*}{Theorem}
\newtheorem*{thma*}{Theorem A}
\newtheorem*{thmaa*}{Theorem A'}
\newtheorem*{thmb*}{Theorem B}
\newtheorem*{thmo*}{Theorem 1.1}
\newtheorem*{thmc*}{Theorem C}
\newtheorem*{thmd*}{Theorem D}
\newtheorem*{thmf*}{Theorem 4.1}
\newtheorem*{remark*}{Remark}
\newtheorem*{conjecture*}{Conjecture}
\newtheorem*{prop*}{Proposition}
\newtheorem*{lem*}{Basic Lemma}
\newtheorem{thm}{Theorem}[section]

\newtheorem{lem}[thm]{Lemma}
\newtheorem{prop}[thm]{Proposition}

\newtheorem*{proofc*}{Proof of Theorem C}

\def\x{\mathbf{x}}

\def\SL{\rm{SL}}

\begin{document}

\author{Youssef Lazar}
\email{ylazar77@gmail.com}
\date{}

\title{Geometric considerations around the Littlewood conjecture}

\maketitle

\begin{abstract} In this paper we adopt a geometric point of view regarding  a famous conjecture due to Littlewood in diophantine approximation of real numbers. Following the spirit of the geometric theory of continued fractions, we give a sufficient condition for the conjecture to hold. An advantage of our method is that it is effective, in the sense that we know how to construct the eventual solution. 
  \end{abstract}


\section{Introduction}
Around 1930, J.E.  Littlewood conjectured that given any two distinct real numbers $\alpha$ and $ \beta$,  the sequence $(n\| n \alpha \| \| n \beta \|)_{n \geq 0}$ can take arbitrary small values. Here $ \| . \|  $  will always denote the function distance  to the nearest integer,  $\| x \| =d(x, \mathbb{Z})$. In the litterature this conjecture is often stated in the following form,
  
 $$~~~~~~ \liminf_{n\rightarrow \infty} ~~n \| n \alpha \|  \| n \beta \| =0.$$

\noindent Despite of its apparent simplicity, this problem remains still open and we are in the situation where there is no explicitely known example, except obvious cases, of a single pair $(\alpha, \beta)$ for which the conjecture is true. Nevertheless, according to some recent results, it is reasonable to think that the conjecture has a very high chance to be true. A first evidence is the discovery of a class of solutions to the problem found by Cassels and Swinnerton-Dyer who proved that the conjecture is true for any pair of irrational cubic numbers lying in the same cubic field \cite{CS55}. This problem sank into oblivion until Polington and Velani proved recently that given a fixed badly approximable number $ \alpha $,  one could find a subset of real numbers $\mathbf{Bad}(\alpha)$ of full Hausdorff dimension such that,  for any element $\beta $ of $\mathbf{Bad}(\alpha) $ the Littlewood is satisfied for the pair $(\alpha, \beta ) $ \cite{PV00}. In the wake of this work,  a tremendous amount of activity have surrounded the problem in the framwork of the metric theory of numbers. New classes of pairs for which the conjecture is valid have been discovered by  Bugeaud, Adamecewski \cite{ABug06} and de Mathan  \cite{DM03} mainly  based on properties of continued fractions.  Bearing some analogy with another conjecture due to Oppenheim and proved by G.A. Margulis, the Littlewood conjecture is a consequence of a conjecture made by Margulis in the context of homogenous dynamics of the space of lattices. The dynamical argument was already implicitely present in the work of Cassels and Swinnerton-Dyer dealing with orbits instead of measures.  Which is considered as the best known result towards the conjecture is due to  Einsiedler, Katok and Lindenstrauss  who solved a conditional version of the Margulis conjecture. As a consequence,  they were able to show that  the set of pairs $(\alpha, \beta ) $ which do not satisfy the Littlewood conjecture has Haussdorff dimension zero \cite{EKL} which is a very strong support to the conjecture. It is worth mentionning the unexpected success of the homegenous dynamics in solving long standing problems in number theory and especially those coming from diophantine approximation.  Another fertile ground was the treatement of the metric diophantine approximation using measure theory and fractal geometry. Understanding the rate of approximation of badly approximable numbers by rationals is the main obstruction to prove the conjecture. For the more details about the state of the art we refer the reader to the recent survey  \cite{BRV16} and also  \cite{Q09}.  Concerning the dynamical point of view see e.g.  \cite{Venk08}.
\subsection*{The main result}
 The aim of this paper is to translate the conjecture into a lattice point problem. Let us be given two badly approximable numbers $(\alpha, \beta)$, the conjecture is then equivalent to the following statement 
 $$ m(f): = \inf_{(x,y,z) \in \mathbb{Z}^{3}} \{  |f(x,y,z)| , x \neq 0 \} =0$$
 where $f(x,y,z)=x(\alpha x -y)(\beta x -z)$. In analytic terms, this assertion equivalent the following one,  for any arbitrary $\varepsilon >0$ there exists $(x,y,z) \in \mathbb{Z}^{3}$, $x\neq 0$ such that 
 $$0< |f(x,y,z)| \leq \varepsilon.$$
  
Let us fixed a positive real number $ \varepsilon$ arbitrarily small.  One is reduced to find a lattice point in the domain bounded by the two cubic level sets $\{ f=\pm \varepsilon\}$. The  very important remark is that this domain, whatever small $\varepsilon$ is,  always contains the intersection of the planes $\{ z=\beta x\}$ and $\{ y=\alpha x\}$ which is nothing else than the line $\mathbb{R}(1,\alpha, \beta)$. Obviously this line cannot afford a lattice point since $1,\alpha, \beta$ are linearly independent over $\mathbb{Q}$. However, given a integer $N$ greater  than $(2\varepsilon)^{-1}$ , a simultaneous version of Dirichlet's Theorem implies the existence of a lattice point $P_{0}= (x_{0}, y_{0}, z_{0}) \in \mathbb{N}^{3}$ such that $1 \leqslant n \leqslant N$ and 
  $$\left\lbrace \begin{array}{c}
\displaystyle |\alpha x_{0}-y_{0}| \leq1/\sqrt{N} \leq \sqrt{2\varepsilon} \\
 \displaystyle |\beta x_{0}-z_{0}| \leq 1/\sqrt{N}\leq \sqrt{2\varepsilon}.
\end{array}\right.$$

\noindent The lattice point $P_{0}$ is then close to the line $\mathbb{R}(1,\alpha, \beta)$, however it not enough close in order to be inside $0< |f| \leq \varepsilon.$. The reason for that is because $\alpha$ and $\beta$ are in $\mathbf{Bad}$. Let us consider the line $P_{0} + \mathbb{R}(1,\alpha, \beta)$ passing through 
$P_{0}$ for the same reason as before this line has no chance to contain a lattice point and let us consider the rational approximation of this line  
  $$(L_{\alpha, \beta}^{n})\left\lbrace \begin{array}{ccc}
x_{n}(t)&=  & x_{0}-t \\
 y_{n}(t)&=  & y_{0}- t c_{2n}(\alpha)\\
 z_{n}(t)&=  & z_{0}- t c_{2n}(\beta)
\end{array}\right.$$
 where $c_{2n}(\alpha)= p_{2n}(\alpha)/q_{2n}(\alpha)$ (resp. $c_{2n}(\beta)$) are the convergents of order $2n$ of $\alpha$ (resp. $\beta$). The reason of the choice of even indices is to take advantage that the positivity of the error term $e_{2n}(\alpha )=\beta - c_{2n}(\alpha)$ (resp. $e_{2n}(\beta )=\alpha - c_{2n}(\beta)$). The line $(L_{\alpha, \beta}^{n}) $  is rational and contains a particular lattice point which corresponds to the time $ t_{n}= \mathrm{lcm}(q_{2n}(\alpha), q_{2n}(\beta))$. The question is whether this lattice point is contained in $\{0< |f| \leq \varepsilon\}$, or in other words,  if
 
 $$0< |f(x_{n}(t_{n}),y_{n}(t_{n}),z_{n}(t_{n}))| \leq \varepsilon.$$ 
 
\noindent  This could be achieved if the line $L_{\alpha, \beta}^{n}$ spends in  $\{0< |f| \leq \varepsilon\}$ a sufficient amount of time in order to reach this lattice point. The most important quantity in our strategy is the time of first entry of the line in  $\{0< |f| \leq \varepsilon\}$, it can be computed as a root of a cubic equation with variable $t$. This leads to very length computations and which will be the subject of a future work. In order to faciltate our task we need to deal with intersection of lines with quadrics for instance lying inside $\{0< |f| \leq \varepsilon\}$. It is not obvious how to do it, but we are fortunate to have a possible candidate given by the following half-cone,

 $$C(N,\varepsilon)= \left\lbrace (x,y,z)\in  [1,N] \times \mathbb{R}^2   ~|~    (\alpha x-y)^2 + (\beta x -z)^2 \leq 2 \varepsilon \dfrac{(N-x)^2}{N(N-1)^2} \right\rbrace.$$   
 
 \noindent The crucial fact that $ C(N,\varepsilon) $ is included in $\{0< |f| \leq \varepsilon\}$ (Proposition \ref{inclusion}) and so we are reduced to solving a quadratic equation instead of a cubic one. This choice induces another problem which is that the cone has as a basis at $x=1$,  a circle of radius $\dfrac{2\varepsilon}{N}$, which tends to be smaller as $N$ goes larger and thus it can happen that the line $L_{\alpha, \beta}^{n}$ never cross the cone $C(N, \varepsilon)$.   To avoid this problem we need to introduce a transversality condition in Proposition \ref{existencetau}, under this assumption,   the time $\tau_{n}$ of first entry of the line $L_{\alpha, \beta}^{n}$ in the cone, is well-defined. Our main result gives a sufficient condition in order to have a lattice point on $L^{n}_{\alpha, \beta}$ after its entry in the cone and in turn it gives a sufficent condition for Littlewoos's conjecture to hold. The main idea is that since both $\alpha$ and $\beta$ are badly approximable numbers, their relative denominators of the converges grows exponentiallly (meaning here not too fast) so that we can bound $t_{n}=\mathrm{lcm}(q_{2n}(\alpha), q_{2n}(\beta))$ in terms of 
 $$ 2^{n-1} \leqslant t_{n} \leqslant  \lambda^{2n}$$
 where $\lambda$ is a common upper bound on the partial quotients of $\alpha$ and $\beta$. Using all those facts, we are able to prove the following theorem 
 
   \begin{thm*}\label{} Let $\varepsilon>0$ be an arbitrary small real number and  $(\alpha,\beta) \in \mathbf{Bad}^{2}$.  If there exists some $n=n(\varepsilon)$ depending on $ \varepsilon $, such that $$ N\left( \sqrt{N}  - \dfrac{1}{\sqrt{N}}\right) \leq \sqrt{2\varepsilon}.
\left( \max\{ e_{2n}(\alpha),  e_{2n}(\beta) \}\right)^{-1}$$ and    $$\tau_{n}\leqslant 2^{n-1} < \lambda^{2n} \leqslant x_{0}-2$$    
where $\tau_{n}$ is the entry time of $L^{n}_{\alpha,\beta}$ in $C(N,\varepsilon)$. Then the Littlewood conjecture holds for $(\alpha, \beta)$.
\end{thm*}

\noindent \textbf{Remarks.}  A clear advantage of this method is its effectiveness, indeed we know how to locate the lattice point solution, namely  $\gamma_{n}(t_{n})$,  if such $n$ exists. However, we are still in the situation where we cannot even produce a single pair which satisfies the hypothesis of the theorem and consequently the Littlewood conjecture. To illustrate this deep difficulty and the limitations of our result, we give in the last section an example of pairs of "very" badly approximable numbers with partial quotients beeing less or equal than three e.g. $\sqrt{2}$, $\sqrt{3}$, $\phi$ the golden ratio etc... For this class we show that it is not possible to ensure simultaneously the two condtions of the theorem.  The main reason seems to be that the cone does not cover enough space inside $\{0 < |f| \leqslant \varepsilon\}$, since its basis tends smaller when $N$ gets larger. The next step to explore consists in trying to replace $\tau_{n}$ the entry time of the line $L^{n}_{\alpha,\beta}$ in the cone by the entry time in the whole $\{0< |f| \leq \varepsilon\}$. There are two clear reasons for doing this choice, the first one is that we can discard the tranversality condition since we are sure that the line will cross $\{0< |f| \leq \varepsilon\}$. The next good point is the better control for the entry time which will be necessarily smaller that our $\tau_{n}$. To see that, it suffices to look to the figure (\ref{default}). We have tried to avoid this by introducing a cone inside $\{0 < |f| \leqslant \varepsilon\}$ and as a product we have obtained this sufficient condition. The hope is that if one considers the life of the line $L^{n}_{\alpha,\beta}$ in  $\{0< |f| \leq \varepsilon\}$  it might be possible to discard the transversality conditions and thus to prove the conjecture. 

\section{Geometric interpretation of the Littlewood conjecture}
 
 \subsection{Lattice points problems for cubic forms associated to the  Littlewood conjecture }In this section we assume that $\alpha$ and $\beta$ are two irrational number such that $1,\alpha, \beta$ are linearly independent over $\mathbb{Q}$. Moreover one can also assume that both $\alpha$ and $\beta$  are badly approximable numbers lying in the unit interval $(0,1) $. As we have seen the validity of the Littlewood conjecture is equivalent to the statement that for every real number $ \varepsilon > 0$, there exists $(x,y,z) \in \mathbb{N}\times \mathbb{Z}^{2}$ such that 
 \begin{equation}\label{conj1}
   |f(x,y,z)|= |x(\alpha x - y)(\beta x -z)| \leq \varepsilon.
\end{equation}
 Geometrically this amounts to find a lattice point $(x,y,z)$ with $x$ positive in the subset of $\mathbb{R}^{3}$ enclosed between the two cubic hypersurfaces $\{f=0 \}$ and $\{f=\varepsilon \}$.  It is clear that since $\alpha$ and $\beta$ are both irrational that $f(x,y,z) \neq 0$ for any $(x,y,z) \in \mathbb{N}\times \mathbb{Z}^{2}$. Let us denote $\{\mid f \mid \leq \varepsilon \} = \cup_{-\varepsilon \leq a \leq \varepsilon} \{ f= a \} $,  one has to prove that 
 $$ \{f \leq \varepsilon \} \cap ( \mathbb{N}\times \mathbb{Z}^{2}) \neq \emptyset.$$
 We consider for any pair of real number $(\alpha, \beta )$  the product of linear forms  $f(x,y,z) = x(\alpha x - y ) (\beta x - z)$ and we fix a pair of integers $(y_0 , z_0 ) \in \mathbb{N}$.\\
Let us define the one variable polynomial of degree 3 $$P(x) =(\alpha \beta)^{-1} f(x,y_0 ,z_0) = x ( x - \alpha^{-1} y_0 ) (x -\beta^{-1} z_0) \in \mathbb{R}[x]$$  The Littlewood conjecture holds for $(\alpha, \beta)$ if for every $\varepsilon >0$ we can find a integer $x_{\varepsilon} \in \mathbb{N}$ such that $|f(x_{\varepsilon})| \leq  \varepsilon$. In order words, we are reduced to find a positive integer in the locus $S_{\varepsilon} := \{x \in \mathbb{R} :  |f(x)| \leq \varepsilon\}$. The roots of $P$ are given by $x_1=0, x_2=  \alpha^{-1} y_0 $ and $x_{3} = \beta^{-1} z_0$. Here we assume that $\alpha, \beta$ are irrational numbers otherwise the result becomes trivial, in particular $x_{1}$ and $x_{2}$ are also irrational numbers. The study of the set of the small values of one variable polynomials is well developped,   a seminal result which is due to H. Cartan, known in the litterature as Cartan's estimate tells us that there exist a covering of $S_{\varepsilon}$ consisting in at most three intervals $I_{1},I_{2}, I_{3} $ respectivily with length $l_{1}, l_{2}, l_{3}$ so that their sums is equal to $e \varepsilon^{1/3}$. In particular, we have 
$$ \lambda (S_{\varepsilon}) \leq \lambda(I_1) +  \lambda(I_2) +  \lambda(I_3) \leq  2e \varepsilon^{1/3}.$$
Hence for $ \varepsilon \leq 1/8e^{3} $, we necessarily have that the length $l_{1}, l_{2}$ and $l_3$ are less that $1$ and therefore each interval $I_j$ for $j=1,2,3$ contains at most one integer and $S_{\varepsilon}$ contains in the best case only 3 integers.  Note that since $0 \in S_{\varepsilon}$, we get  $|S_{\varepsilon} \cap \mathbb{N}| \leq 2$. From this argument we can infer the following crude upper for the number of integral solutions of $|f(x,y,z)| \leq \varepsilon$, 

$$|\{ (x,y,z) \in \mathbb{N} \times [-N,N]^2 : |f(x,y,z)| \leq \varepsilon  \}|  \ll 2 N^2.$$

This is an illustration of the absence of uniform distribution of integral points for regions $\{0< |f| \leqslant \varepsilon\}$. This contrasts with the successful case where the cubic form $f$ is replaced by indefnite quadratic forms in the Oppenheim conjecture.
\noindent \subsection{Parametrization of the level sets of the Littlewood cubic }
 Given a real number $0 < a \leq \varepsilon$, the level set $\{ f=a\}$ has the following parametric formulation

$$\{ f=a\}  ~~~:~~ z = \beta x - \dfrac{a}{x(\alpha x-y)}.$$

\noindent Showing that there exists at least one lattice point on the surface $\{ f=a\}$ seems very complicated using this parametrization. It is not irreducible it splits into several connected components, we are going to work in the first octant, we can assume that the levels are connected.  We are going to makes a change of basis by considering the classical unimodular matrix 
$$ m_{\alpha, \beta} = \left( \begin{array}{ccc}
1 & 0 & 0 \\
\alpha &-1&0 \\
\beta & 0 & -1
\end{array} \right).$$
Given any $(x,y,z) \in \mathbb{R}^{3}$, one has $$m_{\alpha, \beta} (x,y,z)^{t} = (x, \alpha x-y, \beta x -z).$$
Thus if we define the polynomial function $f_{0}$ by setting $f_{0}(x,y,z)=xyz$, we get the fundamental relation $$ f(x,y,z) = f_{0}(m_{\alpha, \beta} (x,y,z)^{t} ).$$
In particular one has the equality 
$$ f(\mathbb{Z}^{3}) = f_{0}(\Lambda_{\alpha, \beta})$$
where $\Lambda_{\alpha, \beta}$ is the unimodular lattice $\mathbb{Z}\left( \begin{array}{c}
1 \\
\alpha\\
\beta
\end{array}\right) \oplus \mathbb{Z}^{2}$ which can been seen as an element of ${\SL}_{3}(\mathbb{R})/{\SL}_{3}(\mathbb{Z})$. The diophantine problem in (\ref{conj1}) is equivalent to the following statement 
\begin{equation}\label{conj2}
 \{0<| f_{0}| \leq \varepsilon \} \cap \Lambda_{\alpha, \beta} \neq \emptyset. 
\end{equation}
In other words one has to find a lattice vector $(x,y,z) \in \Lambda_{\alpha, \beta}$ with  $0 <xyz \leq \varepsilon$.  The symmetry of the situation allows us to reduce ourselves to restrict the research of a lattice vector  to the following domain $$\mathcal{D}(\varepsilon) = \{(x,y,z) \in \mathbb{R}_{+}^{3} ~|~ xyz \leq \varepsilon, \}.$$
Thus we are reduced to prove that for every $\varepsilon >0$, $$\mathcal{D}(\varepsilon) \cap   \Lambda_{\alpha, \beta} \neq \emptyset. $$

\noindent The latter problem does not provide any substantive improvement and we are going to shift to another counting problem.  Studying the distribution of lattice points in convex bodies, paralleotopes for instance is quite well-understood but in general this problem is extremelly difficult. Indeed the domain  $ \mathcal{D}(\varepsilon)$ is not a convex body and one is led to recast the problem in other terms.

\subsection{Elliptic cones inside $ \{|f| \leq \varepsilon \}$} An interesting idea is to introduce a subset of $\mathcal{D}(\varepsilon)$ for which the lattice point counting shows easier. At the first sight, it seems absurd that a domain of smaller volume will be able to produce a lattice point in an easier way. The point is that counting problems in domains bounded by quadric surfaces like hyperboloids or cones have been extensively treated in the litterature and open more perspectives. Of course, such subset must be as large as possible inside $\mathcal{D}(\varepsilon)$. Let us fix  $\varepsilon >0$, and for every positive integer $N>1$ let us define the following region 
$$\mathcal{C}_{0}(\varepsilon, N) : = \{(x,y,z) \in [1,N] \times [0,\sqrt{2\varepsilon}]^2   ~|~    \dfrac{y^2}{2\varepsilon} + \dfrac{z^2}{2\varepsilon} \leq \dfrac{(N-x)^2}{N(N-1)^2} \}.$$
Geometrically $\mathcal{C}_{0}(\varepsilon, N)$ is the region bounded by the positive part of the cone of vertex $(N,0,0)$ and with base given by the circle of equation $y^2 + z^2 =2 \varepsilon/N$ on the plane $x=1$. 



\noindent The following proposition which is at the core of our strategy, shows that this cone is inside the region $\{0 <|f_{0}| \leq \varepsilon\} \cap \mathbb{R}_{+}$, only the inclusion will be useful. 
\begin{lem}\label{inclusion1} For every positive integer $N>1$, one has $$\mathcal{C}_{0}(\varepsilon, N) \subset \mathcal{D}(\varepsilon)$$
Moreover, the half cone  $\mathcal{C}_{0}(\varepsilon, N)$ is not tangent to $ \mathcal{D}(\varepsilon)$.
\end{lem}

\noindent \textbf{Proof.}  The inclusion quite immediate, indeed suppose that $(x,y,z) \in \mathcal{C}_{0}(\varepsilon, N)$, and write
$$xyz = \frac{1}{2}x \left(  y^{2} + z^{2} - (y-z)^{2}  \right).$$

\noindent Therefore
$$xyz \leq \frac{1}{2}x \left(  y^{2} + z^{2}  \right) \leq \dfrac{x(N-x)^2}{N(N-1)^2} \varepsilon .$$
Using that $1 \leqslant x \leqslant N$, we obtain the first assertion.

\noindent Now let us consider the intersection of  $\mathcal{D}(\varepsilon)$ with the plane $x=1$ on which lies the base of our cone, the corresponding equation is given by 
$$0 \leq  z \leq \dfrac{\varepsilon}{y}.$$
The cicular base of our cone has to be tangent to the graph of the hyperbola $\displaystyle z = \dfrac{\varepsilon}{y}$, thus its radius $r$ must satisfies the equations
\begin{center}
$r=\sqrt{y^2 + z^2 }  $ and  $\displaystyle z=\dfrac{\varepsilon}{y}$.
\end{center}
Putting together we get the equation $ r = \sqrt{y^2 +  \left(\dfrac{\varepsilon}{y}\right)^2 } $ and thus after squaring we  are led to 
$$ y^4 -r^2 y + \varepsilon^{2} =0.$$
Set $Y=y^2$, thus $ Y^2 -r^2 Y + \varepsilon^{2} =0.$ A circle centered at origin in the plane may cut the curve $yz=\varepsilon$ at either no points, two points or one point in the first octant. The tangency condition imposes that the quadratic equation must have a unique solution, which reads $\Delta = r^{4}-4 \varepsilon^2 =0$. In other words, the unique solution correspond to the choice $r=\sqrt{2\varepsilon}$.
Thus the base of the required cone in the plane $x=1$, is given by the equation $y^{2}+ z^{2} =2\varepsilon$. Thus this disc of equation $y^2 + z^2 \leq 2 \varepsilon$ is bounded by the curve $yz=\varepsilon$ and the y-axis. Thus, 
$$\mathcal{C}_{0}(\varepsilon, N) \cap \{x=1\} \subset \mathcal{D}(\varepsilon) \cap \{x=1\}.$$
The coordinates of the point of tangency is given by solving $Y^2=-\dfrac{r^2}{2} $, and therefore $y= \sqrt{\varepsilon}$ and also $z=\sqrt{\varepsilon}$. Thus the coordinates of the point of tangency are $(1,\sqrt{\varepsilon}, \sqrt{\varepsilon})$.\\
Now let us consider the cone directed by the circle $y^2 + z^{2} =2\varepsilon$ at $x=1$ and with vertex $(N,0,0)$. Thus the equation of this cone of radius $\sqrt{2\varepsilon}$ and heigth $N-1$ is given by $$  \dfrac{y^2}{2\varepsilon} + \dfrac{z^2}{2\varepsilon} \leq \dfrac{(N-x)^2}{N(N-1)^2}.$$ 

\noindent The equation of the line generating the cone in the plane $y=0$ is the one which passes through the points $(x,z)=(1,\sqrt{\varepsilon})$ and $(x,z)=(N,0)$. The slope of this line equals $- \displaystyle \dfrac{\sqrt{\varepsilon/N}}{N -1}$ and thus the line generating the cone at $y=0$ has an equation which is given by  $$z =- \dfrac{\sqrt{\varepsilon/N}}{N -1} x + \dfrac{\sqrt{\varepsilon/N}N}{N -1}.$$

\noindent Now it remains to show that the line is always below the curve $z=\dfrac{\varepsilon}{xy}$, let us prove that the difference $\Delta z$ is always positive for $1 \leq x \leq N$.  One has,

$$\Delta z = \dfrac{\varepsilon}{xy} + \dfrac{\sqrt{\varepsilon/N}}{N -1}(x-N)   $$

\noindent Since the line is decrasing it suffices to check it for $y<\sqrt{\varepsilon}$ its maximal value. therefore 
$$ \Delta z > \sqrt{\varepsilon/N}(1 + \dfrac{x}{N-1}-\dfrac{N}{N-1}) >\sqrt{\dfrac{\varepsilon}{N}}\dfrac{(x-1)}{N-1} \geq 0.$$

 \noindent This finishes the proof of the Lemma. 
\begin{flushright}
$\square$
\end{flushright}
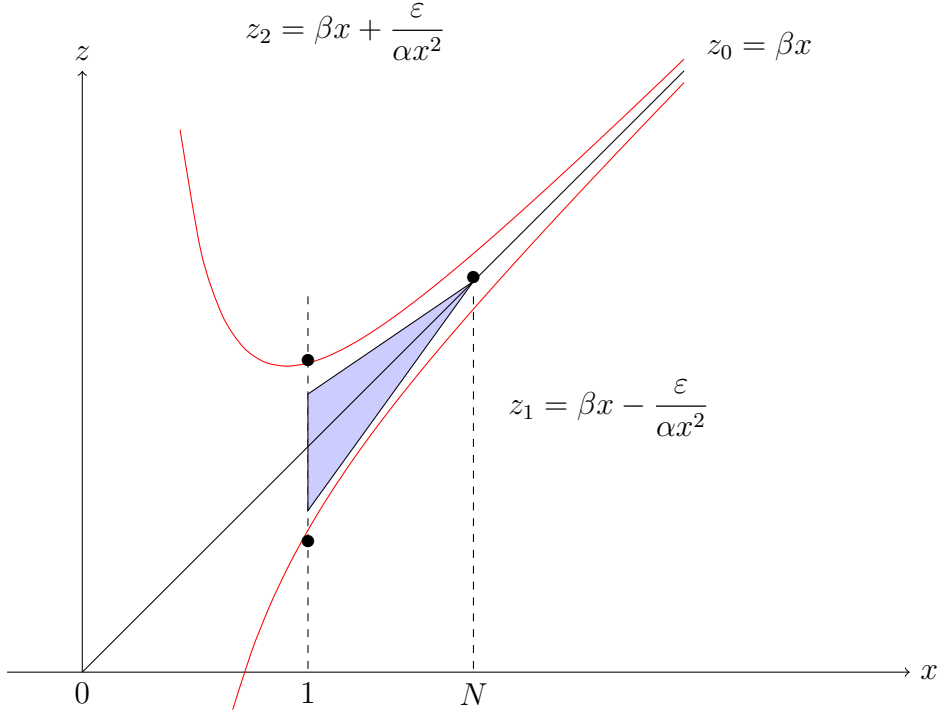
\begin{figure}[htbp]
\begin{center}

\begin{tikzpicture}
\draw[->] (-1,0) -- (11,0);
\draw (11,0) node[right] {$x$};
\draw [->] (0,0) -- (0,8);
\draw (0,8) node[above] {$z$};

\draw (0,0) node[below] {$0$};
\draw[domain=2:8,scale=1, red] plot [smooth](\x,\x -10/\x^2);
\draw[domain=1.3:8,scale=1, red] plot[smooth] (\x,\x +10/\x^2);

\draw [dashed] (3,5) -- (3,0) node[below] {$1$};
\draw [dashed] (5.2,5) -- (5.2,0) node[below] {$N$};

\draw  (7,3) node[above] {$\displaystyle z_{1}= \beta x - \dfrac{\varepsilon}{\alpha x^2}$};
\draw  (3.5,8) node[above] {$\displaystyle z_{2}= \beta x + \dfrac{\varepsilon}{\alpha x^2}$};

\draw  (9,8) node[above] {$\displaystyle z_{0}= \beta x                                                                                   $};

\draw[dotted] (3,1.5) node[above] {$\bullet$};
\draw[dotted] (3,3.9) node[above] {$\bullet$};
\draw[dotted] (5.2,5) node[above] {$\bullet$};

\filldraw[draw=black,fill=blue!20] (3,3.7) -- (5.2,5.2) --(3,2.15)-- cycle;


\draw[domain=0:8,scale=1] plot (\x,{\x});

\end{tikzpicture}
\caption{In blue the section of the cone $ C(N,\varepsilon) \subset \{|f| \leq \varepsilon\}$ in the $xz$-plane.}
\label{default}
\end{center}
\end{figure}

 \noindent A direct consequence of Lemma \ref{inclusion1} is that it supplies a lower bound for the number of points of $\mathcal{D}(\varepsilon) \cap \Lambda_{\alpha, \beta}$, that is, 
$$ \left|\mathcal{C}_{0}(\varepsilon, N)  \cap \Lambda_{\alpha, \beta}\right| \leq \left| \mathcal{D}(\varepsilon) \cap \Lambda_{\alpha, \beta}\right|.$$
The latter inequality is true for every positive integer $N$, and if one is able to find an $N$ such that $ \left| \mathcal{C}_{0}(\varepsilon, N)  \cap \Lambda_{\alpha, \beta} \right| > 0$ then this is sufficient to validate the conjecture. As it is, the counting problem in the lattice  $\Lambda_{\alpha, \beta}$ is seems also very hard to apprehend.  To move to a counting in the standard lattice we use the following linear transformation 
$$  \left|\mathcal{C}_{0}(\varepsilon, N)  \cap \Lambda_{\alpha, \beta}\right| = \left|  m_{\alpha, \beta}^{-1} \mathcal{C}_{0}(\varepsilon, N)  \cap  ( \mathbb{N}\times \mathbb{Z}^{2} )\right| = \left|  m_{\alpha, \beta} \mathcal{C}_{0}(\varepsilon, N)  \cap  ( \mathbb{N}\times \mathbb{Z}^{2} )\right|$$
the last equality holds  since $m_{\alpha, \beta}^{-1}=m_{\alpha, \beta}$. We are reduced to finding a lattice point in the domain $ m_{\alpha, \beta}^{-1} \mathcal{C}_{0}(\varepsilon, N) $ which we are going to make explicit
$$ m_{\alpha, \beta}^{-1} \mathcal{C}_{0}(\varepsilon, N) = \{ (x,y,z)\in   [1,N] \times [0,\sqrt{2\varepsilon}]^2  ~|~    m_{\alpha, \beta}(x,y,z)^{t} \in \mathcal{C}_{0}(\varepsilon, N)     \}. $$

\noindent Thus $ m_{\alpha, \beta}^{-1} \mathcal{C}_{0}(\varepsilon, N) $ is the is equal to $$C(N,\varepsilon)= \left\lbrace (x,y,z)\in  [1,N] \times \mathbb{R}^2   ~|~    (\alpha x-y)^2 + (\beta x -z)^2 \leq 2 \varepsilon \dfrac{(N-x)^2}{N(N-1)^2} \right\rbrace.$$  
Since 
Thus we have shown the following crucial inclusion,

\begin{prop}\label{inclusion} For every positive integer $N>1$, one has $$ C(N,\varepsilon) \subset  \{0 <|f| < \varepsilon \}. $$
\end{prop}
\vspace{1cm}
\noindent \textbf{Remarks.} Hence a sufficient condition for conjecture (\ref{conj1}) to hold is the existence of some $N$ and  $(x,y,z) \in ([1,N]\cap \mathbb{N}) \times \mathbb{Z}^2$ such that $(x,y,z) \in  C(N,\varepsilon).$  
\noindent Indeed the latter would ensure that $$ 0 <|m_{\alpha, \beta}^{-1} \mathcal{C}_{0}(\varepsilon, N) \cap  ( \mathbb{N}\times \mathbb{Z}^{2} ) | \leq |\{ 0 < |f|  \leq \varepsilon\} \cap  ( \mathbb{N}\times \mathbb{Z}^{2} )|. $$

\noindent Geometrically $C(N,\varepsilon)$ represents the interior of  a cone with circular base and with its axis directed by the vector $(1, \alpha, \beta)$. The vertex of the cone $C(N,\varepsilon)$ is $(N, N \alpha, N \beta)$. We obviously have $\dfrac{(N-x)^2}{N(N-1)^2} \leq 1/N$ for every $1 \leqslant x \leqslant N$,  the subset $C(N,\varepsilon) $ lies inside the parralelepiped, 
$$P(N,\varepsilon)= \{ (x,y,z)\in  [1,N] \times \mathbb{R}^2  ~|~    |\alpha x-y|,  |\beta x-y| \leq \sqrt{2\varepsilon}/N  \}. $$\\ The subset  $P(N,\varepsilon) $ constrains any point $(x,y,z)$ to be closed to the both the planes $y=\alpha x $ and $z= \beta x$ or in other words to their intersection which is exactly the line $\mathbb{R}(1,\alpha, \beta)$. The end of the paper will be dedicated to the study of the sequence  $$N \mapsto |m_{\alpha, \beta}^{-1} \mathcal{C}_{0}(\varepsilon, N) \cap  ( \mathbb{N}\times \mathbb{Z}^{2} ) | .$$ Therefore the Littlewood problem (\ref{conj1}) reduces to finding an integer $ N(\varepsilon) $ at which this sequence assumes a positive value.


\section{Lattice points on lines of best approximation}
\subsection{Best approximation of badly approximable numbers by continued fractions}

 \noindent Let us recall some of the main properties of continued fractions we are going to need, see for instance Chapter 3 \cite{EW} and \cite{Sch80} for more details and also \cite{OK} for the geometric theory.
 
 The usual representation of  a real number $\alpha$  as a continued fraction takes the following form

 $$\alpha= [a_0;a_1,a_2,a_3, \ldots]= \displaystyle a_0+ \dfrac{1}{ a_1+ \dfrac{1}{ a_2+ \dfrac{1}{a_3+ ...}} } .$$

\noindent where $(a_{k})$ is the sequence positive integers which characterizes $\alpha$, and called the partial quotients of $\alpha$. In the same way,  we denote by $\beta = [b_0;b_{1},b_2,b_{3}, \ldots]$ the corresponding continued fraction expansion. Since we are assuming that $\alpha  $ and $ \beta $ are irrational numbers, their expansion is infinite. The convergent of order $n$ of $ \alpha $ is given by the following fraction
$$  c_{n}(\alpha) =[a_0;a_1,a_2,a_3, \ldots, a_{n}].$$

\noindent Let us denote by $ \dfrac{p_{n}(\alpha)}{q_{n}(\alpha)} $ the reduced expression of  $c_{n}(\alpha)$, the terminology is justified by  the fact that $\lim_{n} c_{n}(\alpha)= \alpha$. The integers $p_{n}$ and $q_{n}$ enjoy the following recursive property 
\begin{equation}
p_{n}(\alpha) = a_{n}p_{n-1}(\alpha) + p_{n-2}(\alpha)
\end{equation}
\begin{equation}\label{reccurence-qn}
q_{n}(\alpha) = a_{n}q_{n-1}(\alpha) + q_{n-2}(\alpha).
\end{equation}

\noindent with initial values $ p_{-1}=1 $, $ q_{-1} =0$, $ p_{0} =a_{0}$ and $ q_{0} =1$. 
Also one has the following relation $p_{n}(\alpha) q_{n-1}(\alpha)- p_{n-1}(\alpha) q_{n}(\alpha) = (-1)^{n+1}$ which implies the convergents forms an alternating sequence in the following sense 
\begin{equation}\label{evenconvergents}
c_{0} < c_{2} <  \ldots < c_{2n}< \ldots < \alpha < \ldots , c_{2n+1}< \ldots < c_{3} < c_{1}. 
\end{equation}

\noindent The error in the approximation of $\alpha$ by $c_{n}(\alpha)$ is denoted by $e_{n}(\alpha)$, i.e. $e_{n}(\alpha)= \alpha - c_{n}(\alpha)$, it is exactly given by 
$$ e_{n}(\alpha) =   \dfrac{(-1)^{n}}{\alpha_{n+1} q_{n}(\alpha)+q_{n+1}(\alpha)}$$

\noindent where $\alpha_{n+1}$ is such that $\alpha =[a_{0}; a_{1}, \ldots, a_{n}, \alpha_{n+1}] $ (see Lemma 3 E, \cite{Sch80}).
It is not difficult to see that $\alpha_{n+1} = a_{n+1} + 1/\alpha_{n+2} \leq a_{n+1} +1$. From that we deduce  the following classical inequalities ( \cite{EW} ex. 3.1.5 p. 76 )
\begin{equation}\label{error1}
 \dfrac{1}{2q_{n}(\alpha)q_{n+1}(\alpha)} \leq    |e_{n}(\alpha)| \leq  \dfrac{1}{q_{n}(\alpha)q_{n+1}(\alpha)} \leq \dfrac{1}{q_{n}(\alpha)^{2}}
\end{equation}

\noindent The property (\ref{evenconvergents}) ensures that both $e_{2n}(\alpha)$ and $e_{2n}(\beta)$ are positive and  (\ref{error1}) reads
\begin{equation}\label{error22}
 \dfrac{1}{2q_{2n}(\alpha)q_{2n+1}(\alpha)} \leq    e_{2n}(\alpha) \leq  \dfrac{1}{q_{2n}(\alpha)q_{2n+1}(\alpha)} \leq \dfrac{1}{q_{2n}(\alpha)^{2}}.
\end{equation}
\subsubsection*{Growth of the denominators of the convergents, Levy's Theorem}
The set of badly approximable numbers can be defined briefly as the set, 

$$\mathbf{Bad} := \{ \alpha \in \mathbb{R} ~\vert~ \inf_{q \geqslant 1} q\|q\alpha \| > 0 ~\mathrm{for}~\mathrm{every}~q\in \mathbb{N}    \}.$$
This numbers corresponds to the real numbers for which the Dirichlet theorem cannot be improved up to a mutiplicative constant. More precisely, $\alpha \in \mathbf{Bad}$, if there exists a constant $C$ for every $p$ and $q$ positive integers
$$ \dfrac{C}{q^{2}} < \left| \alpha - \dfrac{p}{q}\right| \leqslant \dfrac{1}{q^{2}}. $$  

The inequality of the left hand side is satisfied for every quadratic irrational number after Liouiville's Theorem, thus those are badly approxiable numbers. It is conjectured that these are the only ones, and it still an open to know whether there exists an algebraic number of degree $>2$ which is in  $\mathbf{Bad}$.  In the other inequality, Roth's  Theorem tells us that this inequality cannot be improved, for every real $\delta >$ there are only a finite number of pairs $(p,q)$ such that 
  $$\left| \alpha - \dfrac{p}{q}\right| \leqslant \dfrac{1}{q^{2+\delta}}.$$
  
\noindent   Roughly one can say that if $\alpha$ is badly approximable then as $q$ gets large,
  $$ \|q\alpha \| \asymp \dfrac{1}{q}$$ 
   The notation $f(q) \asymp g(q)$ means that the ratio of $f$ by $g$ is bounded above and below by constants which does not depends on $q$. 
  From the recursive formula (\ref{reccurence-qn}) above it is not difficult to deduce that $\left(( q_{n}(\alpha)\right)_{n\geq 1}$ and  $\left(( q_{n}(\alpha)\right)_{n\geq 1}$ are increasing sequences. Indeed since $\alpha, \beta \in \mathbf{Bad}$, their partial quotient are both bounded so that we can set for $n \geq 1$   $$M :=\max(\sup_{n} a_{n}(\alpha); \sup_{n} a_{n}(\beta)) < \infty$$
\noindent We get for both $\alpha$ and $\beta$, the inequalities

$$q_{n-1} \leqslant q_{n} \leqslant \left( M+1 \right) q_{n-1}.$$

\noindent Hence, we obtain the following bounds 
$$2^{(n-2)/2} \leqslant  q_{n} \leqslant \left( M+1 \right)^{n}.$$
\noindent In particular, setting $\lambda := ( M+1 )^{2}$ 
\begin{equation}\label{lw-qn}
  2^{(n-2)/2} \leq q_{n} \leq \lambda^{n/2}.
\end{equation}
and taking even indices we get 
$$   2^{n-1}  \leqslant  q_{2n}  \leqslant \lambda^{n}.$$
The inequalities in (\ref{error1}) and (\ref{lw-qn}) give 
\begin{equation}\label{error2}
 \dfrac{1}{2\lambda^{2n+1/2}}  \leq \dfrac{1}{2q_{2n}(\alpha)q_{2n+1}(\alpha)}  \leq e_{2n}(\alpha) \leq  \dfrac{1}{q_{2n}(\alpha)q_{2n+1}(\alpha)} \leq \dfrac{1}{q_{2n}(\alpha)^{2}} \leq \dfrac{1}{2^{2n-2}}
\end{equation} 
and the same holds for $\beta$. If $\alpha$ is an irrational number in $(0,1)$,  in order to quantify the rate of growth of the denominators of the convergents, one defines the Levy upper and lower constant as follows
\begin{equation}\label{levy}
\kappa^{\ast}(\alpha) =\limsup_{n} \dfrac{\log q_{n}(\alpha)}{n}  ~~\mathrm{and} ~~~~ \kappa_{\ast}(\alpha) =\liminf_{n} \dfrac{\log q_{n}(\alpha)}{n}.
\end{equation}

In the case when this two quantities coincide we call it the Levy constant of $\alpha $ and we just denotes it $ \kappa(\alpha) $.  Regarding the typical behaviour of the denominators $ q_{n}(\alpha)$, P. Levy proved the following surprizing result that for almost every $\alpha$, 
$$\kappa(\alpha) = \exp \left(\frac{\pi^{2}}{12\log 2} \right).$$
This constant number is called the Levy constant, so Levy's theorem tells us the typical growth of any irrational $a$ number is 
$$ q_{n}(a) \sim  \exp \left(\frac{n\pi^{2}}{12\log 2} \right).$$

\subsection{Time spent by best approximation lines in $ C(N,\varepsilon) $ }
\subsubsection*{Parametrization of lines passing through a given lattice vector}
A simultaneous version of Dirichlet's Theorem tells us that for every $\varepsilon >0$ and for every $N>(2\varepsilon)^{-1}$ there exists a lattice vector $P_{0}=(x_{0},y_{0},z_{0})$ with $1 \leq x_{0} \leq N$

$$\left\lbrace \begin{array}{c}
\displaystyle |\alpha x_{0}-y_{0}| \leq1/\sqrt{N} \leq \sqrt{2\varepsilon} \\
 \displaystyle |\beta x_{0}-z_{0}| \leq 1/\sqrt{N}\leq \sqrt{2\varepsilon}.
\end{array}\right.$$

\noindent This produces a lattice which is near to be as we wish, indeed
$$|f(x_{0},y_{0},z_{0})| =|x_{0}||\alpha x_{0}-y_{0}|  |\beta x_{0}-z_{0}| \leq 2 \varepsilon x_{0} \leq 2\varepsilon N.$$

\noindent The upper bound is $2\varepsilon N$ is greater than one so $ |f(x_{0},y_{0},z_{0})| $ is small only if $x_{0}(\varepsilon)$ is small enough. Since we are dealing with badly approximable numbers there exists a constant $ C $ depending only on $ \alpha $ and $ \beta $ such that 

$$\left\lbrace \begin{array}{c}
\displaystyle C < x_{0} |\alpha x_{0}-y_{0}| \leq \sqrt{2\varepsilon} \\
\displaystyle C < x_{0} |\beta x_{0}-z_{0}| \leq \sqrt{2\varepsilon}.
\end{array}\right.$$
In particular this gives us that 
\begin{equation}\label{badly}
\frac{C}{\sqrt{2\varepsilon}} < x_{0}(\varepsilon) \leq N.
\end{equation}

\noindent There exists a unique line passing through $P_{0}$ and in the direction of the vector $(1,\alpha,\beta)$ we denote it $L_{\alpha, \beta}$. This line is parralel with the axis of the cone  for which we recall the definition 
$$C(N,\varepsilon)= \left\lbrace (x,y,z)\in  [1,N] \times \mathbb{R}^2   ~|~    (\alpha x-y)^2 + (\beta x -z)^2 \leq 2 \varepsilon \dfrac{(N-x)^2}{N(N-1)^2} \right\rbrace.$$
 The parametrization of the part of the line going downwards is given by 
$$L_{\alpha, \beta} = P_{0} + \mathbb{R}_{-}\left(1,\alpha,\beta   \right)=\left\{ (x_{0}-t, y_{0}-t \alpha, z_{0} -t \beta ) \mid  t \geq 0 \right \}$$

\noindent This line cannot contain any other lattice point than $P_{0}$, so we are going to consider a rational approximation of $ L_{\alpha, \beta} $ using a line passing through $ P_{0} $ with direction vector $(1,c_{n}(\alpha),c_{n}(\beta))  $, 
$$L^{n}_{\alpha, \beta} = \left\{ (x_{0}-t, y_{0}-t c_{2n}(\alpha), z_{0} -t c_{2n}(\beta) ) \mid 0 \leq  t \leq x_{0}-1 \right \}$$
A parametrization of the line segment $L^{n}_{\alpha, \beta}$  is given by $ \gamma_{n}(t)=(x_{n}(t), y_{n}(t), z_{n}(t))~~~ (0 \leq t \leq x_{0}-1)$ where 
$$\left\lbrace \begin{array}{ccc}
x_{n}(t)&=  & x_{0}-t \\
 y_{n}(t)&=  & y_{0}- t c_{2n}(\alpha)\\
 z_{n}(t)&=  & z_{0}- t c_{2n}(\beta)
\end{array}\right.$$
The $ L^{n}_{\alpha, \beta}  $ converges to $ L_{\alpha, \beta} $ as $n$ tends to infinity and we interested in the intersection points $(P_{n})_{n}$ of the lines $ L^{n}_{\alpha, \beta}  $ with the boundary of $ C(N,\varepsilon) $. More than that, and roughly speaking we need to estimate the time spent by a mobile point starting from $P_{0}$ moving on the line $ L^{n}_{\alpha, \beta}  $  until it comes across $ A(N,\varepsilon) $ (see figure \ref{dessin2}). This time denoted $\tau_{n}$ is the time of first entry in  $ C(N,\varepsilon) $ and the intersection point will be $ \gamma_{n}(\tau_{n})$. This time will play a central role in what follows and its existence is ensured by the following proposition.

\begin{figure}[htbp]
\label{dessin2}
\begin{center}

\begin{tikzpicture}
\draw[->] (-1,0) -- (11,0);
\draw (11,0) node[right] {$x$};
\draw [->] (0,0) -- (0,8);
\draw (0,8) node[above] {$z$};

\draw (0,0) node[below] {$0$};
\draw[domain=2:8,scale=1, red] plot [smooth](\x,\x -10/\x^2);
\draw[domain=1.3:8,scale=1, red] plot[smooth] (\x,\x +10/\x^2);

\draw [dashed] (3,5) -- (3,0) node[below] {$1$};
\draw [dashed] (5.2,5) -- (5.2,0) node[below] {$N$};
\draw [dashed] (4.5,3.8) -- (4.5,0) node[below] {$x_{0}$};

\draw  (7,3) node[above] {$\displaystyle z= \beta x - \dfrac{\varepsilon}{\alpha x^2}$};
\draw  (3.5,8) node[above] {$\displaystyle z= \beta x + \dfrac{\varepsilon}{\alpha x^2}$};

\draw  (9,8) node[above] {$\displaystyle z= \beta x                                                                                   $};

\draw  (9,7) node[above] {$(L^{n}_{\alpha, \beta})$};

\draw[dotted] (3,1.5) node[above] {$\bullet$};

\draw[dotted] (3,1.5) node[above] {$\bullet$};
\draw[dotted] (3,3.9) node[above] {$\bullet$};
\draw[dotted] (4.7,3.5) node[above] {$\bullet P_{0}$};
\draw[dotted] (4,2.4) node[above] {$\bullet \gamma_{n}(\tau_{n})$};
\draw[dotted] (6.3,4.8) node[above] {$\bullet (N,\alpha N, \beta N)$};

\filldraw[draw=black,fill=blue!20] (3,3.7) -- (5.2,5.2) --(3,2.15)-- cycle;


\draw[domain=0:8,scale=1] plot (\x,{\x});
\draw[domain=0.5:8,scale=1,green] plot (\x,{\x -0.7});


\end{tikzpicture}
\caption{The green line $(L^{n}_{\alpha, \beta})$ passes through $P_{0}$ and cuts the cone $C(N,\varepsilon)$ at $\gamma_{n}(\tau_{n})$ is the tranversality condition is satisfied.}
\label{default}
\end{center}
\end{figure}
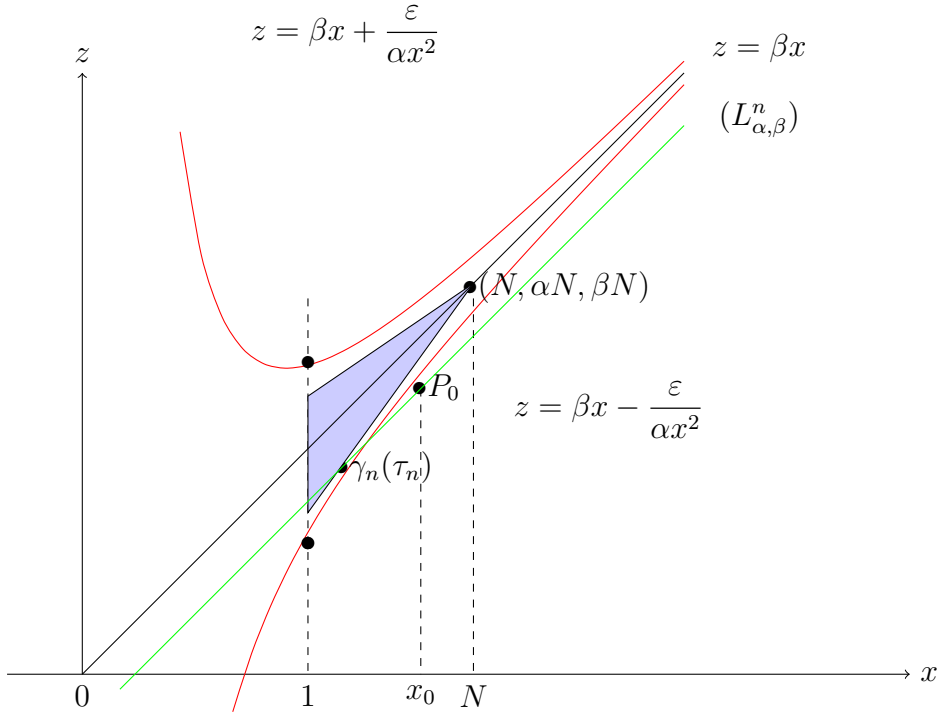
\begin{prop} \label{existencetau}  If $N$ and $ n $ are two integers satisfiing $$ N\left( \sqrt{N}  - \dfrac{1}{\sqrt{N}}\right) \leq \sqrt{2\varepsilon}
\left(2 \max\{ e_{2n}(\alpha),  e_{2n}(\beta) \}\right)^{-1}$$
 then we can find a real constant $0< \tau_{n} < x_{0}$ depending on $N,\varepsilon,n,
\alpha, \beta $ (and $ P_{0} $) such that
\begin{center}
$ \gamma_{n}(t) \in L^{n}_{\alpha, \beta} \cap A(N,\varepsilon)$ if and only if $t \in [  \tau_{n}, x_{0}-1]. $
\end{center}
\end{prop}
\noindent \textbf{Proof.} Asking that $\gamma_{n}(t) \in L^{n}_{\alpha, \beta}\cap C(N,\varepsilon)$ is equivalent to the condition
 $$ \left( \alpha x_{n}(t)-y_{n}(t)\right)^{2} + \left( \beta x_{n}(t)- z_{n}(t)\right)^{2}  \leq \dfrac{2 \varepsilon}{N(N-1)^2}(N-x_{n}(t))^{2}.$$
\noindent  Thus, 
$$  \left( \alpha (x_{0}-t)-(y_{0}-tc_{2n}(\alpha)\right)^{2} + \left( \beta (x_{0}-t)-(z_{0}-tc_{2n}(\beta)\right)^{2} \leq \dfrac{2 \varepsilon}{N(N-1)^2}(N-x_{0}+t)^{2}.$$

\noindent Rearranging the terms we get, 
$$  \left( \alpha x_{0}-y_{0}- t(\alpha -c_{2n}(\alpha))\right)^{2} + \left( \beta x_{0}-z_{0}- t(\beta -c_{2n}(\beta))\right)^{2}  \leq \dfrac{2 \varepsilon}{N (N-1)^2}(N-x_{0}+t)^{2}.$$

\noindent In order to clarify let us introduce some notations. We set $U_0 =  \alpha x_{0}-y_{0}$ and  $V_0 =  \beta x_{0}-y_{0}$ and we put $\displaystyle \varphi(\varepsilon,N) = \dfrac{2 \varepsilon}{N(N-1)^2}$.  Thus, for all $0 \leq t \leq x_{0}-1$

$$ \left( U_{0}- t e_{2n}(\alpha)\right)^{2} + \left( V_{0}-t e_{2n}(\beta)\right)^{2}  \leq \varphi(\varepsilon,N)(N-x_{0}+t)^{2}.$$

\noindent Expanding we are left to study the following quadratic inequation
$$\left\lbrace \varphi(\varepsilon,N) -e_{2n}(\alpha)^{2} - e_{2n}(\beta)^{2}\right\rbrace t^{2} +  2t \left\lbrace \varphi(\varepsilon,N)(N-x_{0}) +e_{2n}(\alpha)U_{0} + e_{2n}(\beta) V_{0} \right\rbrace $$ $$+\left\lbrace \varphi(\varepsilon,N)(N-x_{0})^{2} - U_{0}^{2}- V_{0}^{2} \right\rbrace \geq 0.$$
\noindent The cutting times of the line $L^{n}_{\alpha, \beta}$ with the whole cone $C(N,\varepsilon) = A(N,\varepsilon) \cup \sigma (A(N,\varepsilon))$ where $\sigma$ is the central symmetry with center the vertex $(N,N\alpha, N\beta)$ are exactly the roots in $t$ of the left hand side of the inequality. The discriminant of the quadratic polynomial is given by
$$ D_{n}(N,\varepsilon) =4  \left\lbrace \varphi(\varepsilon,N)(N-x_{0}) +e_{2n}(\alpha)U_{0} + e_{2n}(\beta) V_{0} \right\rbrace^{2}$$  $$- 4\left\lbrace \varphi(\varepsilon,N) -e_{2n}(\alpha)^{2} - e_{2n}(\beta)^{2}\right\rbrace  \left\lbrace \varphi(\varepsilon,N)(N-x_{0})^{2} - U_{0}^{2}- V_{0}^{2}  \right\rbrace. $$
\begin{lem}\label{Dn positive}
For every integer $n>1$, if the value of $N$ satsfies the condition 
$$ N\left( \sqrt{N}  - \dfrac{1}{\sqrt{N}}\right) \leq \sqrt{2\varepsilon}
\left(2 \max\{ e_{2n}(\alpha),  e_{2n}(\beta) \}\right)^{-1}$$
Then  $ D_{n}(N,\varepsilon) >0$, and the line  $L^{n}_{\alpha, \beta}$ cuts $C(N,\varepsilon)$ in a a unique point.

\end{lem}

\noindent\textit{Proof of the Lemma.}  For the sake of brevity we simply write $\varphi$ for $\varphi(N,\varepsilon)$. We  compute the terms,
$$\dfrac{1}{4} D_{n}(\varepsilon, N)= \varphi^{2}(N-x_{0})^{2} + 2\varphi (N-x_{0}) \{e_{2n}(\alpha)U_{0} +e_{2n}(\beta)V_{0} \} $$ $$+ e_{2n}(\alpha)^{2}U_{0}^{2}+ 2e_{2n}(\alpha)e_{2n}(\beta) U_{0} V_{0} + e_{2n}(\beta)^{2}V_{0}^{2} - \varphi^{2}(N-x_{0})^{2}$$
$$+e_{2n}(\alpha)^{2}\varphi (N-x_{0})^{2}+e_{2n}(\beta)^{2}\varphi (N-x_{0})^{2}+\varphi(U_{0}^{2} + V_{0}^{2}) $$
$$ - e_{2n}(\alpha)^{2}(U_{0}^{2} +V_{0}^{2}) - e_{2n}(\beta)^{2}(U_{0}^{2} + V_{0}^{2}).$$

\noindent A first reduction gives, 

$$\dfrac{1}{4} D_{n}(\varepsilon, N) = 2\varphi (N-x_{0})(e_{2n}(\alpha)U_{0} +e_{2n}(\beta)V_{0} ) + 2e_{2n}(\alpha)e_{2n}(\beta) U_{0} V_{0} 
+e_{2n}(\alpha)^{2}\varphi (N-x_{0})^{2}$$ $$+e_{2n}(\beta)^{2}\varphi (N-x_{0})^{2}+\underbrace{\varphi(U_{0}^{2} + V_{0}^{2})  - e_{2n}(\alpha)^{2} V_{0}^{2} - e_{2n}(\beta)^{2}U_{0}^{2}}$$

\noindent For $n$ large enough we have $\varphi -  e_{2n}(\alpha)^{2} \geq 0$ and $\varphi -  e_{2n}(\beta)^{2} \geq 0$, so the expression under brace is nonnegative and takes the following form 
$$  \left( \varphi -  e_{2n}(\alpha)^{2} \right)V_{0}^{2}  +  \left( \varphi -  e_{2n}(\beta)^{2} \right)U_{0}^{2}.$$
We have the following 
\begin{equation}\label{DN1}
\dfrac{1}{4} D_{n}(\varepsilon, N) = 2\varphi (N-x_{0})(e_{2n}(\alpha)U_{0} +e_{2n}(\beta)V_{0} ) + 2e_{2n}(\alpha)e_{2n}(\beta) U_{0} V_{0} 
+e_{2n}(\alpha)^{2}\varphi (N-x_{0})^{2}$$ $$+e_{2n}(\beta)^{2}\varphi (N-x_{0})^{2}+\left( \varphi -  e_{2n}(\alpha)^{2} \right)V_{0}^{2}  +  \left( \varphi -  e_{2n}(\beta)^{2} \right)U_{0}^{2}. 
\end{equation}
 \noindent The latter is also nonnegative, indeed using the assumption 
 $$\sqrt{N} \left(N-1 \right) =N\left( \sqrt{N}  - \dfrac{1}{\sqrt{N}}\right) \leq 
\dfrac{\sqrt{2\varepsilon}}{2\max\{ e_{2n}(\alpha),  e_{2n}(\beta) \}}. $$
 Thus 
 $$ \sqrt{\varphi }=\dfrac{\sqrt{2\varepsilon}}{\sqrt{N} \left(N-1 \right)} \geq 2\max\{ e_{2n}(\alpha),  e_{2n}(\beta) \}$$
 and in particular 
 \begin{center}
 $\varphi -  e_{2n}(\alpha)^{2} \geq 0$ and $\varphi -  e_{2n}(\beta)^{2} \geq 0$.
 \end{center}
 Even better, we have 
 \begin{equation}\label{denominatortau}
 \varphi -  e_{2n}(\alpha)^{2} -  e_{2n}(\beta)^{2} \geq \varphi -2\max\{ e_{2n}(\alpha),  e_{2n}(\beta) \}^{2} \geq  \varphi -4\max\{ e_{2n}(\alpha),  e_{2n}(\beta) \}^{2} \geq 0.
 \end{equation}
 \noindent This shows that $ D_{n}(\varepsilon, N) >0$, and thus line  $L^{n}_{\alpha, \beta}$ cuts the half $C(N,\varepsilon)$ at exactly one point which corresponds to the positive root. This proves the lemma.

\noindent \textit{End of the proof of Proposition \ref{existencetau}}. The cutting times with the entire cone $A(N,\varepsilon) \cup \sigma \left( A(N,\varepsilon)\right)$ where $ \sigma $ is the central symmetry  with respect to the vertex $(N,N\alpha, N\beta)$ are given by by the two distincts roots
\begin{center}
$\displaystyle t_{-} = \dfrac{-2\left( e_{2n}(\alpha)U_{0} + e_{2n}(\beta) V_{0}+ \varphi(\varepsilon,N)(N-x_{0})\right)-\sqrt{D_{n}(N,\varepsilon)}}{ \varphi(\varepsilon,N) -e_{2n}(\alpha)^{2} - e_{2n}(\beta)^{2}  }$
\end{center}
and 
\begin{center}
$\displaystyle t_{+} = \dfrac{-2\left( e_{2n}(\alpha)U_{0} + e_{2n}(\beta) V_{0}+ \varphi(\varepsilon,N)(N-x_{0})\right)+\sqrt{D_{n}(N,\varepsilon)}}{ \varphi(\varepsilon,N) -e_{2n}(\alpha)^{2} - e_{2n}(\beta)^{2} }.$
\end{center}

\noindent \textit{Small angles hypothesis.} ~~ We have already seen that when $n$ is large enough the quantity $D_{n}$ is positive but it can still happen that the line has two intersection points with $A(N,\varepsilon)$. We are going to show that this will never occurs, at least if $n$ is taken large enough. Geometrically, this fact is quite obvious since the line $L_{\alpha,\beta}$ is parralel to $\mathbb{R}(1,\alpha,\beta)$ the axis of the cone 
\noindent Let us define the angle $ \theta_{n}(\alpha, \beta )$ between the lines  $L^{n}_{\alpha, \beta}$ and  $L_{\alpha, \beta}$, it can be computed as 
$$  \theta_{n}(\alpha, \beta )=\cos^{-1} \left(\dfrac{1 + \alpha c_{2n}(\alpha) + \beta c_{2n}(\beta) }{\sqrt{1 + \alpha^{2}  + \beta^{2} }\sqrt{1 +  c_{2n}(\alpha)^{2} +  c_{2n}(\beta)^{2} }}\right).$$

Clearly $\lim_{n}  \theta_{n}(\alpha, \beta ) = 0$, thus it becomes very small as $n$ gets large. In this case, the line  $L^{n}_{\alpha, \beta}$ cuts $C(N,\varepsilon)$ in exactly one point. The  opposite case corresponds to the situation when the line $L_{\alpha, \beta}$  intersects transversally the half-cone $C(N,\varepsilon)$. Since the convergents are assumed to be close to their relative limits, this situation is very unlikely to happen at least for $n>n_{0}$ where $n_{0}$ is large enough. \\
Therefore for $n \geq n_{0}$ we are in the situation when the line  $L^{n}_{\alpha, \beta}$ cuts in exactly one point. This is what we assume by now.\\
\noindent The time $t_{-}$ is the (negative) cutting time of the line with the upper part of the cone symmetric to $C(N,\varepsilon)$  with respect to the vertex $(N,N\alpha, N\beta)$ so as we have seen it has to be discarded. 
Thus $t^{+}$ is the first time of entry of the line $L^{n}_{\alpha, \beta} $ in  the half-cone $C(N,\varepsilon)$ and correspond to the unique point of intersection.  
Hence for $n \geq n_{0}$, we define the quantity
\begin{equation}\label{entrytime}
\tau(\varepsilon, N, n, P_{0}) =t^{+}=  \dfrac{-2\left( e_{2n}(\alpha)U_{0} + e_{2n}(\beta) V_{0}+ \varphi(\varepsilon,N)(N-x_{0})\right)+\sqrt{D_{n}(N,\varepsilon)}}{ \varphi(\varepsilon,N) -e_{2n}(\alpha)^{2} - e_{2n}(\beta)^{2}}.
\end{equation} 
It remains to prove that $t^{+}$ is positive.
\noindent Thus for any $\tau(\varepsilon, N, n, P_{0})\leq t \leq x_{0}-1$, $\gamma_{n}(t) \in C(N,\varepsilon)$ since $1 \leq x_{0} \leq N$. This proves the proposition.
\begin{flushright}
$\square$
\end{flushright}

\subsection{The existence of lattice points in $L^{n}_{\alpha,\beta} \cap A(N,\varepsilon)$}

In the previous section, we used a parametric equation of the portion of the line $L^{n}_{\alpha,\beta} $ lying inside $C(N,\varepsilon)$ and thus also in  $P(N,\varepsilon)$. It is given by $ \gamma_{n}(t)=(x_{n}(t), y_{n}(t), z_{n}(t))~~~ (\tau_{n} \leq t \leq x_{0}-1)$ where 
$$\left\lbrace \begin{array}{ccc}
x(t)&=  & x_{0}-t \\
 y(t)&=  & y_{0}- t c_{2n}(\alpha)\\
 z(t)&=  & z_{0}- t c_{2n}(\beta)
\end{array}\right.$$
where $\tau_{n}$ is as in (\ref{entrytime}). A sufficient condition in order to make this segment containing a lattice point is that the positive integer $t_{n} : = \mathrm{lcm} (q_{2n}(\alpha), q_{2n}(\beta))$ lies in the interval $[\tau_{n}, x_{0}-1]$. Indeed in that case, the following vector will be a lattice point given by 
$$\left\lbrace \begin{array}{ccc}
x_{n}(t_{n})&=  & x_{0}- \mathrm{lcm} (q_{2n}(\alpha), q_{2n}(\beta)). \\ \\
 y_{n}(t_{n})&=  & \displaystyle y_{0}- \dfrac{\mathrm{lcm} (q_{2n}(\alpha), q_{2n}(\beta))}{q_{2n}(\alpha)} p_{2n}(\alpha)\\
  z_{n}(t_{n})&=  & \displaystyle  z_{0}- \dfrac{\mathrm{lcm} (q_{2n}(\alpha), q_{2n}(\beta))}{q_{2n}(\beta)} p_{2n}(\beta).\\
\end{array}\right. $$
Thus if we find such an $n$ then the point $X_{n}:=\gamma_{n}(t_{n})=(x_{n}(t_{n}), y_{n}(t_{n}), z_{n}(t_{n}))$ will be an integral vector in $A(N,\varepsilon)$ with $ x_{n}(t_{n}) \geq 1$. A consequence of that is that one would have  $|A(N, \varepsilon)\cap (\mathbb{N} \times \mathbb{Z}^{2})|>0$ and in particular  $|\{ |f|\leq \varepsilon\}\cap (\mathbb{N} \times \mathbb{Z}^{2})|>0$.  Thus finding such $t_{n}$ is sufficient to establish the truth of the Littlewood conjecture for $(\alpha, \beta)$. \\

\noindent We are reduced to find an integer  $n(\varepsilon)$ depending on $\varepsilon$ such that 
\begin{equation}\label{timespent}
\tau_{n}= \tau(N,n, \varepsilon) \leq  t_{n}= \mathrm{lcm} (q_{2n}(\alpha), q_{2n}(\beta)) < x_{0}.
\end{equation}
Indeed if $t_{n} < \tau_{n}$, the vector $\gamma_{n}(t_{n})=(x(t_{n}), y(t_{n}), z(t_{n})) $ will produce a latiice point in $ P(N,\varepsilon )$ but stiil outside $A(N,\varepsilon)$. In the other hand, if $t_{n} \geq x_{0}$ then  one has $x_{n}({t_{n}}) = x_{0}-t_{n}\leq 0$ and thus $(x_{n}(t_{n}), y_{n}(t_{n}), z_{n}(t_{n})) $ would also live outside  $A(N,\varepsilon)$. Note that the integral vector $P_{0}=(x_{0},y_{0},z_{0}) \in P(N,\varepsilon)$ depends on $\varepsilon$, and a fortioti $x_{0}$ depends on it too.

\subsubsection*{Some Facts}  Let us collect some important facts which are going to be useful thereafter,

\begin{itemize}
\item[$(1)$] We have that $N $ should be at least equal to  $(2\varepsilon)^{-1}$ because of the application the simultaneous version of Dirichlet's theorem supplies an integral vector $P_{0}=(x_{0},y_{0}, z_{0})$ such that
$\displaystyle |U_{0}|= |\alpha x_{0} - y_{0}| \leq \dfrac{1}{N^{1/2}} \leq \sqrt{2\varepsilon}$,    and $\displaystyle  |V_{0}|=|\beta x_{0} -z_{0}| \leq \dfrac{1}{N^{1/2}} \leq \sqrt{2\varepsilon}$.
\item[$(2)$]  Since $ \alpha $ and  $ \beta $  are badly approximable we have that $\dfrac{C}{\sqrt{2\varepsilon}} \leqslant x_{0} \leqslant N$  where  $$ C =\max\{\inf_{q\geq1} q\| q\alpha\|,\inf_{q\geq1}q\| q\beta\|\}>0. $$

\item[$(3)$] \label{fact4} Using the fact that the partial quotients of the two badly approximable numbers $\alpha$ and $\beta$ are bounded, say for both by $M$, one can find a universal constant $ \lambda \geqslant 2 $ such that  $$   2^{n}  \leqslant t_{n}=\mathrm{lcm} (q_{2n}(\alpha), q_{2n}(\beta)) \leqslant \lambda^{2n}.$$

The integer $ t_{n}=\mathrm{lcm} (q_{2n}(\alpha), q_{2n}(\beta))$ is at least $q_{2n}(\alpha)$ or $q_{2n}(\beta)$ if both are equal, thus $2^{n-1} \leqslant t_{n}$. In the other hand, $t_{n}$ is at most the product $q_{2n}(\alpha) q_{2n}(\beta)$, and this bound is reached as soon as $\mathrm{gcd} (q_{2n}(\alpha), q_{2n}(\beta)) =1$. To sum up
$$2^{n-1} \leqslant t_{n} \leqslant \lambda^{2n}.$$

\noindent \item[$(4)$] \label{fact5}The sequence of integers $(t_{n})_{n \geqslant 1}$ grows exponentially with the order of convergence $n$.  Indeed since the sequence $\left( \dfrac{\log t_{n}}{n} \right)_{n}$ is bounded one can find a positive real number $\delta=\delta(\alpha,\beta)$  such that
$$  \lim_{n} \dfrac{\log t_{n}}{n} =\log  \delta.$$  
\noindent Therefore, we have the following estimate as $ n $ gets large,
\begin{equation}
t_{n} \sim \delta^{n}.
\end{equation}
\end{itemize}

\vspace{0.5cm}
\subsection*{A sufficent condition for the Littlewood conjecture to hold }
Now we reach the key result showing that if for some $ n $ large enough the line $L^{n}_{\alpha,\beta}  $ spends a sufficient amount of time in $C(N,\varepsilon)$ then there is a lattice point in $C(N,\varepsilon)$ which is equivalent to the Littlewood conjecture.

\begin{thm}\label{thmsuff} Let $\varepsilon>0$ be an arbitrary small real number and  $(\alpha,\beta) \in \mathbf{Bad}^{2}$.  If there exists some $n=n(\varepsilon)$ depending on $ \varepsilon $, such that $$ N\left( \sqrt{N}  - \dfrac{1}{\sqrt{N}}\right) \leq \sqrt{2\varepsilon}.
\left( \max\{ e_{2n}(\alpha),  e_{2n}(\beta) \}\right)^{-1}$$ and    $$\tau_{n}\leqslant 2^{n-1} < \lambda^{2n} \leqslant x_{0}-2$$    
where $\tau_{n}$ is the entry time of $L^{n}_{\alpha,\beta}$ in $C(N,\varepsilon)$ as in (\ref{entrytime}). Then the Littlewood conjecture holds for $(\alpha, \beta)$.
\end{thm}
\noindent \textit{Proof of the Theorem.} Suppose that such an $n$ exists, then  because of Proposition \ref{existencetau} the transversality condition implies that the entry time of $L^{n}_{\alpha,\beta}$ in $C(N,\varepsilon)$ is well defined. The second condition and the fact (3) gives us the following 
$$   \tau_{n}\leqslant 2^{n-1} \leqslant  t_{n} \leqslant \lambda^{2n} \leqslant x_{0}-2.$$
where $t_{n}= \mathrm{lcm}(q_{2n}(\alpha), q_{2n}(\alpha))$.
Thus the lattice point $$\gamma_{n}(t_{n}) = ( x_{0}- t_{n},   \displaystyle y_{0}- \dfrac{t_{n}}{q_{2n}(\alpha)} p_{2n}(\alpha),  z_{0}- \dfrac{t_{n}}{q_{2n}(\beta)} p_{2n}(\beta) )$$
lies in $L^{n}_{\alpha,\beta} \cap C(N,\varepsilon)$. In particular, using proposition \ref{inclusion} we obtain that  $$0 < |f(\gamma_{n}(t_{n}))| \leqslant \varepsilon$$
which just says that the Littlewood conjecture holds for $(\alpha,\beta)$.
\begin{flushright}
$\square$
\end{flushright}
\section{The difficulty of finding a lattice point in $C(N,\varepsilon)$}

For conveniance we will restrict to pairs of badly approximable vector whose partial quotients consists only of $1$ and $2$, and $3$ e.g. the golden ratio $\phi$, $\sqrt{2}, \sqrt{3}$...  we denote  $(\alpha,\beta) \in \mathbf{B}(3)^{2}$ where 
$$\mathbf{B}(N)= \{ \alpha \in  \mathbb{R}  ~|~ \alpha= [a_{0}; a_{1}, a_{2}, a_{3}, \ldots ] , ~ 1 \leqslant a_{j} \leq N, j=0,1, \ldots \}. $$ We are going to see that it is difficult to find a lattice point in the cone $C(N,\varepsilon)$ corresponding to $(\alpha, \beta)$. Let $ \varepsilon >0 $, we are want to prove that there exists a lattice vector $(x,y,z)\in \mathbb{N} \times \mathbb{Z}^{2}$ such that 
$$ 0 < |f(x,y,z)|= |x(\alpha x -y)(\beta x -z)| \leq \varepsilon. $$
After Proposition \ref{inclusion}  it would be sufficient to finding a lattice point in $ C(N, \varepsilon) $ for some $N \geq1$. Let us take $ N>(2\varepsilon)^{-1} +1$, then the simultaneous Dirichlet's approximation theorem applied to $(\alpha,\beta)$ gives us an integral vector $P_{0}=(x_{0},y_{0},z_{0})$ such that 
$$\left\lbrace \begin{array}{c}
1 \leqslant x_{0}  \leqslant N \\
 |\alpha x_{0} - y_{0}| \leq  \sqrt{2\varepsilon}\\
  |\beta x_{0} - z_{0}| \leq  \sqrt{2\varepsilon}.
 \end{array}\right. $$
 
\noindent Let us consider $P_{0}\left(\varepsilon \right)=(x_{0}, y_{0},z_{0}) \in P(N,\varepsilon) \cap (\mathbb{N} \times \mathbb{Z}^{2})$ is the lattice vector produced by the Dirichlet Theorem depending on the choice of $N$. By Theorem \ref{thmsuff}  it suffices to prove the existence of an integer $n(\varepsilon)$ such that
$$\tau_{n}= \tau(N,n, \varepsilon)  \leq 2^{n-1}.$$

\noindent We need a precise upper estimate for
 $$ \tau_{n}=  \dfrac{-2\left( e_{2n}(\alpha)U_{0} + e_{2n}(\beta) V_{0}+ \varphi(\varepsilon,N)(N-x_{0})\right)+\sqrt{D_{n}(N,\varepsilon)}}{ \varphi(\varepsilon,N) -e_{2n}(\alpha)^{2} - e_{2n}(\beta)^{2}}. $$
 
\noindent  \noindent The assumption on $ \alpha, \beta$ tells us that we should have 
$$M :=\max(\sup_{n} a_{n}(\alpha); \sup_{n} a_{n}(\beta)) = 3.$$

\noindent In particular $\lambda =(M+1)^2 = 16=2^4$, after replacing we get 

\begin{equation}\label{errorsqauare}
 \dfrac{1}{\lambda^{4n+4}} = \dfrac{1}{2^{16(n+1)}}   \leq e_{2n}(\alpha)^{2}+e_{2n}(\beta)^{2} \leq \dfrac{1}{2^{4n-5}}
\end{equation}

\noindent and 

 \begin{equation}\label{errorzero}
\dfrac{C}{\lambda^{2n+1/2}x_{0}} = \dfrac{C}{2^{8n+4}x_{0}}  \leq e_{2n}(\alpha)U_{0} + e_{2n}(\beta)V_{0}\leq \dfrac{2 \sqrt{2\varepsilon}}{2^{2n-2}}.
\end{equation}

\noindent Let us focus on the denominator of $  \tau_{n} $, from the transversality condition we already know that
$$ \varphi =\dfrac{2\varepsilon}{N \left(N-1 \right)^{2}} \geq 4\max\{ e_{2n}(\alpha),  e_{2n}(\beta) \}^{2}. $$

\noindent Using (\ref{error2}) again, we infer that
$$\varphi- e_{2n}(\alpha)^{2}-e_{2n}(\beta)^{2}  \geq 4\max\{ e_{2n}(\alpha),  e_{2n}(\beta) \}^{2} -e_{2n}(\alpha)^{2} - e_{2n}(\beta)^{2} \geq \dfrac{4}{(2\lambda^{2n+1})^{2}}  -   \dfrac{2}{(2\lambda^{2n+1})^{2}}   .$$

\noindent Thus, 
$$\varphi- e_{2n}(\alpha)^{2}-e_{2n}(\beta)^{2}  \geq \dfrac{1}{2\lambda^{4n+4}}.$$

\begin{equation}\label{denotau2}
\left( \varphi- e_{2n}(\alpha)^{2}-e_{2n}(\beta)^{2} \right)^{-1}   \leq 2^{16n+17}.
\end{equation}

\noindent Now we need an upper estimate of $D_{n}^{1/2}$,
\begin{lem}\label{Dn estimate}
Let us set $X=2\varepsilon$ and $ u=2^{n} $ then one has the estimate

$$ D_{n}(\varepsilon, N)^{1/2} \ll   4\sqrt{2} X u^{-2}+\sqrt{2}~  X^{5/2} + \sqrt{2}  X^{4}  u^{2}/4+2 \sqrt{2} X u^{2}.$$
\end{lem}
\noindent \textit{Proof of the Lemma.}   It remains to estimate $D(n,\varepsilon)$ in this range, using (\ref{DN1}) we begin with

\begin{equation}\label{DN2}
\dfrac{1}{4} D_{n}(\varepsilon, N) = 2\varphi (N-x_{0})\left\lbrace e_{2n}(\alpha)U_{0} +e_{2n}(\beta)V_{0} \right\rbrace + 2e_{2n}(\alpha)e_{2n}(\beta) U_{0} V_{0} $$ $$+\varphi (N-x_{0})^{2} \left\lbrace e_{2n}(\alpha)^{2}+e_{2n}(\beta)^{2}\right\rbrace +\left( \varphi -  e_{2n}(\alpha)^{2} \right)V_{0}^{2}  +  \left( \varphi -  e_{2n}(\beta)^{2} \right)U_{0}^{2}. 
\end{equation}

\noindent Also using  (\ref{lw-qn})  we get $$e_{2n}(\alpha)U_{0} + e_{2n}(\beta)V_{0}  \leq \sqrt{2\varepsilon}\left( \dfrac{1}{q_{2n}(\alpha)^{2}} +   \dfrac{1}{q_{2n}(\beta)^{2}}\right)  \leq \dfrac{2 \sqrt{2\varepsilon}}{2^{2n-2}}$$
and 
$$e_{2n}(\alpha)^{2} + e_{2n}(\beta)^{2}  \leq \left( \dfrac{1}{q_{2n}(\alpha)^{4}} +   \dfrac{1}{q_{2n}(\beta)^{4}}\right)  \leq \dfrac{1}{2^{4n-5}}.$$

\noindent From the definition of $P_{0}$ arising from Dirichlet's Theorem with  (\ref{error1}) and  (\ref{error2})  we get  $$\left( \varphi -  e_{2n}(\alpha)^{2} \right)V_{0}^{2}  +  \left( \varphi -  e_{2n}(\beta)^{2} \right)U_{0}^{2} \leq 4\varepsilon \left( \varphi - \dfrac{1}{2\lambda^{2n+1}}  \right).$$  
\noindent In this range one has, 
 $$ D_{n}(\varepsilon, N) \leq \dfrac{4\sqrt{2 \varepsilon}\varphi (N-x_{0})}{2^{2n-2}}+ \dfrac{\varphi (N-x_{0})^{2}}{2^{4n-5}}+ 4\varepsilon \left(\varphi +\dfrac{1}{2^{4n-4}} -  \dfrac{1}{2\lambda^{2n+1}}   \right).$$
 \begin{flushright}
$\square$
\end{flushright}

\noindent Now,  using the estimate for $D_{n}(\varepsilon, N)$ obtained in Lemma \ref{Dn estimate}, and replacing $ \lambda$ we infer that
 
 $$ D_{n}(\varepsilon, N) \leq \left( 2^{6n+7}    4\sqrt{2 \varepsilon}\varphi (N-x_{0}) + \varphi (N-x_{0})^{2}2^{4n+10} +4\varepsilon \varphi 2^{8n+5}+ 4\varepsilon 2^{8n+9}-4\varepsilon \right) 2^{-8n-5}. $$
 
 $$ D_{n}(\varepsilon, N) \leq \left( 2^{6n+9}   \sqrt{2 \varepsilon}\varphi (N-x_{0}) + \varphi (N-x_{0})^{2}2^{4n+10} +\varepsilon \varphi 2^{8n+7}+ \varepsilon 2^{8n+11} \right) 2^{-8n-5}. $$
\noindent Dividing $2^{4n+5}$, we get 
  $$ D_{n}(\varepsilon, N) \leq \dfrac{ 2^{2n+4}   \sqrt{2 \varepsilon}\varphi (N-x_{0}) + 32 \varphi (N-x_{0})^{2} +\varepsilon \varphi 2^{4n+2}+ \varepsilon 2^{4n+6} }{2^{4n}}. $$
\noindent   Using the fact that $N > (2\varepsilon)^{-1}+1$   (Dirichlet's condition)
 \begin{equation}\label{fi}
 \varphi (N-x_{0})  \leq  \dfrac{2\varepsilon}{N(N-1)} \leq (2\varepsilon)^{3} ~~~\mathrm{and}~~~ \varphi (N-x_{0})^{2}  \leq  \dfrac{2\varepsilon}{N} \leq (2\varepsilon)^{2}
\end{equation}

  $$ D_{n}(\varepsilon, N) \leq \dfrac{ 2^{2n+4}  (2 \varepsilon)^{7/2} + 32(2\varepsilon)^{2} +(2\varepsilon)^{5}  2^{4n+2}+ (2\varepsilon) 2^{4n+5} }{2^{4n}}. $$
  
\noindent   Dividing by $2^{5}$ 
  
  $$ D_{n}(\varepsilon, N) \leq (2\varepsilon)^{2} \dfrac{ 1+2^{2n-1}  (2 \varepsilon)^{3/2} +(2\varepsilon)^{3}  2^{4n-3}+ (2\varepsilon)^{-1} 2^{4n} }{2^{4n-5}}. $$ 
  Using $\sqrt{1+x} < 1+ \dfrac{1}{2}x$ in the numerator, we obtain 
  
$$ D^{1/2}_{n}(\varepsilon, N) \leq (2\varepsilon)  \dfrac{ 1+2^{2n-2}  (2 \varepsilon)^{3/2} +(2\varepsilon)^{3}  2^{4n-4}+ (2\varepsilon)^{-1} 2^{4n-1} }{2^{2n-5/2}}. $$ 

\noindent Let us set for convenience $X=2\varepsilon$ and $u=2^{n}$, thus
$$ D^{1/2}_{n}(\varepsilon, N) \leq X \dfrac{ 1+u^{2}  X^{3/2}/4 +X^{3}  u^{4}/16+ X^{-1} u^{4}/2 }{u^{2}/4\sqrt{2}}. $$ 
\noindent  Therefore we arrive to the following bound for the square root of  $D_{n}$
   $$ D^{1/2}_{n}(\varepsilon, N) \leq  4\sqrt{2} X u^{-2}+\sqrt{2}~  X^{5/2} + \sqrt{2}  X^{4}  u^{2}/4+2 \sqrt{2} X u^{2}. $$ 
  
\begin{flushright}
$  \square$
\end{flushright}

\noindent We come back to $\tau_{n}$;  remembering that the application of Dirichlet's requires that $N>(2\varepsilon)^{-1}+1$ beside the tranversality condition, we get from the inequalities (15-19) 
$$ \tau_{n} \ll  \left( D^{1/2}_{n}(\varepsilon, N)  - \dfrac{C}{2^{8n+3}x_{0}}- 8\max\{ e_{2n}(\alpha),  e_{2n}(\beta) \}^{2} (N-x_{0})  \right)(\varphi - \dfrac{1}{\lambda^{4n+1}})^{-1}.$$

\noindent Using (\ref{denotau2})
$$ \tau_{n} \ll  \left( D^{1/2}_{n}(\varepsilon, N)  - \dfrac{C}{2^{8n+3}x_{0}}- 8\max\{ e_{2n}(\alpha),  e_{2n}(\beta) \}^{2} (N-x_{0})  \right) 2^{16n+17}.$$

\noindent We have, 
$$8 \max\{ e_{2n}(\alpha),  e_{2n}(\beta) \}^{2} \geq \frac{8}{4 \lambda^{4n+1}} = \dfrac{1}{2^{16n+2}}$$
so 

$$ \tau_{n} \ll  \left( D^{1/2}_{n}(\varepsilon, N)  - \dfrac{C}{2^{8n+3}x_{0}}-  \dfrac{N-x_{0} }{2^{16n+2}}  \right) 2^{16n+17}.$$

\noindent Using the Lemma
\begin{equation} \label{polyinu}
 \tau_{n} \ll  \left( 4\sqrt{2} X u^{-2}+\sqrt{2}~  X^{5/2} + \left( \sqrt{2}  X^{4}/4 +2 \sqrt{2}  \right)  u^{2} - \dfrac{C}{8 u^{8}x_{0}}-  \dfrac{N-x_{0} }{4 u^{16}}  \right) u^{16} 2^{17}.
\end{equation}

\noindent Now the theorem \ref{thmsuff} affirms that the pair $(\alpha, \beta)$ to satisfies Littlewood conjecture, if there exists some $n$ which fullfills the transversality condition such that
\begin{center}
$ \tau_{n}  \leq 2^{n-1}  < 2^{4n} \leq x_{0}-1$ 
\end{center}
or with our notations 
\begin{equation}\label{rangeu}
\tau_{n}  \leq u/2  < u^{4} \leq x_{0}-1.
\end{equation} 

\noindent Let us definie the following polynomial with correspond to the right hand side in (\ref{polyinu})
$$ \psi_{\varepsilon}(u) : =   2^{17+\frac{1}{2}} \left(   X^{4}/4 +2  \right)  u^{18}+ 2^{17+\frac{1}{2}}~  X^{5/2} u^{16} + 4\sqrt{2} X u^{14}- \dfrac{2^{14}C}{x_{0} } u^{8}- 2^{15}(N-x_{0}) .$$

\noindent Thus the inequality (\ref{polyinu}) can be written as
$$   \tau_{n} \ll  \psi_{\varepsilon}(2^{n}).  $$

\noindent Let us first focus on the set of positive real solutions of the inequation
 $$    \psi_{\varepsilon}(u)< \dfrac{u}{2}. $$  
 
\noindent The transversality condition says that $n$ should at least satisfies
 
 \begin{equation}\label{transversal}
 N\left( \sqrt{N}  - \dfrac{1}{\sqrt{N}}\right) \leq \sqrt{2\varepsilon}
\left(2 \max\{ e_{2n}(\alpha),  e_{2n}(\beta) \}\right)^{-1} \leq \sqrt{2\varepsilon} ~2^{8n+4}= 16 \sqrt{2\varepsilon} ~u^{8}.
\end{equation} 
 
 \noindent In view of (\ref{rangeu}) $u$ should be at least in range
 $$\left( \dfrac{\sqrt{N}\left( N-1 \right)}{16\sqrt{2\varepsilon}} \right)^{1/8} \leq u \leq  (x_{0}-1)^{1/4}.$$
 
\noindent   Again the Dirichlet condition  $N > (2\varepsilon)^{-1}+1$   (Dirichlet's condition) gives us the range
 
 $$  \dfrac{1}{2^{5/8} (2\varepsilon)^{3/16}}  \leq u  \leq (x_{0}-1)^{1/4}.$$
 
\noindent  The corresponding range in term of $n$ is 
 
 $$ -\dfrac{\log (2^{5/8} (2\varepsilon)^{3/16})}{\log 2}  \leqslant  n  \leqslant \dfrac{\log (x_{0}-1)}{4\log 2}.$$
 
 \noindent We are reduced to find such integer $n$ such that $ \psi_{\varepsilon}(2^{n}) \leq 2^{n-1} $. It is more convenient to see $\psi_{\varepsilon}$ as a polynomial function with real variables in the range given above, and to study the function $h_{\varepsilon}(u)=\psi_{\varepsilon}(u)-\frac{1}{2}u$. Therefore we may solve the inequation $h_{\varepsilon}(u) \leq 0$, the large degree does not allows to compute exact solutions but we can use some calculus.
 
\noindent Let us write $h_{\varepsilon}(u)= au^{18} + b u^{16} + c u^{14} - d u^{8} -  u/2 -e $ where $a,b,c,d$ and $e$ are the coefficients as in its definition of $\psi_{\varepsilon}$. Thus, solving $h_{\varepsilon}(u) \leq 0$ amounts to solve after dividing by $u^{14}$

\begin{equation}\label{ineqh}
 a u^{4} + bu^{2} + c \leq \dfrac{d}{u^6} +  \dfrac{1}{2u^{13}}+ \dfrac{e}{u^{14}}
\end{equation}. 
 
\noindent Here $d=2^{14}C/x_{0}$ and $ e = 2^{15}(N-x_{0})$, we know from (\ref{transversal}) that $  d=2^{14}C/x_{0} \leq 2^{14} X$ and $$  e \leq 2^{15}( N-1) \leq \dfrac{\sqrt{2\varepsilon}}{2\sqrt{N}} \max\{ e_{2n}(\alpha),  e_{2n}(\beta) \}^{-1} \leq 8u^{8}. $$ 

\noindent Also using $u^{-1} \leq 2^{5/8} X^{3/16}$, we get an estimate for the right hand side of (\ref{ineqh}) is given by 

$$ \dfrac{2^{14}X}{u^6} +  \dfrac{2^{-1}}{u^{13}}+ \dfrac{8}{u^{6}} \leq 2^{14}X( 2^{5/8} X^{3/16})^{6}  + 2^{-1}( 2^{5/8} X^{3/16})^{13} + 2^{3}( 2^{5/8} X^{3/16})^{6}$$

\noindent We are reduced to study the following system the inequalities 

\begin{equation}\label{ineqh2}
\left\lbrace \begin{array}{c}
  a u^{4} + bu^{2} + c \leq 2^{7/4}X^{1/8} + 2^{57/8} X^{39/16} + 2^{27/4} X^{9/8}.\\
 2^{-5/8} X^{-3/16}  \leq u  \leq (x_{0}-1)^{1/4}.
  \end{array}\right. 
\end{equation}
\noindent Let us set $ Z=u^{2} $, thus we have to solve 

 $$2^{35/2}(X^{4} +2) Z^{2} + 2^{35/2} X^{5/2} Z + (4\sqrt{2}X-2^{7/4}X^{1/8} + 2^{57/8} X^{39/16} + 2^{27/4} X^{9/8}) \leq  0. $$

\noindent Now there is no need to go further, since that 

$$2^{35/2}(X^{4} +2) Z^{2} \leq 2^{7/4}X^{1/8} -2^{35/2} X^{5/2} Z - (4\sqrt{2}X+ 2^{57/8} X^{39/16} + 2^{27/4} X^{9/8}) . $$

\noindent Indeed in the range given by (\ref{ineqh2}), we have that the right hand side of the inequality just above is negative while the left hand side remains always positive, this lead to a contradiction.

\section{Final remarks}  As we have seen, finding a lattice point in a cone shaped domain inside  $\{0 < |f| \leq \varepsilon \}$ using lines has not shown a satisfactory issue. However,  the situation is not completely compromized even if there is a very tiny chance that this cone gives a satisfactory answer to this lattice point problem. Forgetting about cones, there is a real hope that this strategy could succeed if one considers the intersection of $L_{\alpha, \beta}^{n}$ with the whole  $\{0 < |f| \leq \varepsilon \}$ rather than a subset of it. In that case, we will be sure that the entry time will be shorter that the $\tau_{n}$ of the main theorem. As a consequence, the line will spend more time in   $\{0 < |f| \leq \varepsilon \}$ which leaves more chance to bound above it by $t_{n}$. This can be done by using the same method.\\
 Given a real number $\varepsilon >0$, for any $N > (2\varepsilon)^{-1}$ one can consider the lattice point  $P_{0}$ arising from Dirichlet's Theorm and which is close to  the line $ \mathbb{R} (1,\alpha, \beta)$ as in \S 3.1.  Therefore asking that $\gamma_{n}(t) \in L^{n}_{\alpha, \beta}\cap \{0 < |f| \leq \varepsilon \}$ is equivalent to the condition
 $$ -\varepsilon   \leq x_{n}(t) \left( \alpha x_{n}(t)-y_{n}(t)\right) \left( \beta x_{n}(t)- z_{n}(t)\right)  \leq \varepsilon.$$
Taking account of the symmetry, we only focus on the first octant, thus one is reduced to solve for  $1\leq t \leq x_{0}$ the following inequation
 $$   0  < (x_{0}-t) \left( te_{2n}(\alpha) - U_{0}\right) \left( t e_{2n}(\beta)-V_{0}  \right) \leq \varepsilon.$$
The line $L_{\alpha, \beta}^{n}$ intersects the boundary $\{ |f| = \varepsilon \}$ in at most three points, but the choice the line which is nearly parallel to the  asymptote irrational line $\mathbb{R}(1,\alpha, \beta)$ ensures that $L_{\alpha, \beta}^{n}$ cuts $\{ f= - \varepsilon\}$ in a unique point. Here the situation is much more delicate than the quadratic situation. In fact, finding the root of a cubic polynomial in one variable is not a difficult task but the estimate of the first time of entry $$ \inf\{t \in \mathbb{R}_{+} ~:~  \gamma_{n}(t) \in L_{n}^{\alpha, \beta}  \cap \{0 < |f| \leq \varepsilon \}  \}$$
is much involved in this situation.  Indeed, one has to deal with cubic roots corresponding to the cubic equation
$$ (t-x_{0}) \left( te_{2n}(\alpha) - U_{0}\right) \left( t e_{2n}(\beta)-V_{0}  \right) = \varepsilon$$
and estimating the entry time needs more efforts. Even harder evaluations are, the avalaiblity of closed forms makes possible to get sharp estimates of the entry time. This strategy will be detailled in a forthcoming work.

\end{document}